\documentstyle[12pt]{amsart}

\makeatletter                                            
\renewcommand{\@begintheorem}[2]{                        
\rm \trivlist \item [\hskip \labelsep {\bf #2\ \ #1.}]   
                                }                        
\makeatother                                             

\newcommand{\ts}{\vspace{\baselineskip}\noindent{\bf Proof.}$\;\;$}
\newcommand{\ZZ}{{\bf Z}}
\newcommand{\QQ}{{\bf Q}}

\newcommand{\CC}{{\bf C}}

\newcommand{\PP}{{\bf P}}

\newcommand{\bP}{{\Bbb P}}  

\newcommand{\bZ}{{\Bbb Z}}
\newcommand{\cO}{{\cal O}}
\newcommand{\inj}{\hookrightarrow}

\newcommand{\ra}{\rightarrow}

\begin{document}

\title{Half twists and the cohomology of hypersurfaces}
\author{B. van Geemen}
\address{Dipartimento di Matematica, Universit\`a di Pavia,
  via Ferrata 1, I-27100 Pavia, Italia}
\email{geemen@@dragon.ian.pv.cnr.it}
\author{E. Izadi}
\address{Department of Mathematics, Boyd Graduate Studies Research
  Center, University of Georgia, Athens, GA 30602-7403, USA}
\email{izadi@@math.uga.edu}

\begin{abstract}
A Hodge structure $V$ of weight $k$ on which a CM field acts defines, under 
certain conditions, a Hodge structure of weight $k-1$, its half twist.
In this paper we consider hypersurfaces in projective space with a cyclic
automorphism which defines an action of a cyclotomic field on a Hodge 
substructure in the cohomology. We determine when the half twist exists and 
relate it to the geometry and moduli of the hypersurfaces. We use our results 
to prove the existence of a Kuga-Satake correspondance for certain cubic 
4-folds.
\end{abstract}

\maketitle

{\bf Introduction}

\

Given a rational Hodge structure $V$ of weight $k$ with an action of a
CM-field $K$ such that no element of $K$ has complex conjugate
eigenvalues on $V^{k,0}$, one can associate to it a polarized Hodge
structure $V_{1/2}$ (the half twist of $V$) of weight $k-1$ with the
same underlying vector space. This half twist depends on the choice
of a CM-type for $K$ and was defined by one of us in \cite{geemen01}.
It has the property that there is an inclusion of Hodge structures
$V_{1/2}\subset V\otimes H^1(A_K)$ (up to Tate twist) where $A_K$ is
an abelian variety with CM by $K$. The Hodge conjectures translate
into interesting problems on the geometry of
varieties having Hodge structures which allow a half twist, some of
which we consider in this paper.

The half twist is also related to Kuga-Satake varieties. In case $V$ has
weight $2$ with $V^{ 2,0}$ of dimension $1$, Kuga and Satake
\cite{kugasatake67} define an abelian variety $KS(V)$ such that
$H^1(KS(V))$ is the even part of the Clifford algebra of the
polarization of $V$. They show that $V$ is a direct summand of $H^1
(KS(V) )^{\otimes 2}$. If an imaginary quadratic field $K$ acts on
$V$ then $V_{1/2}$ and $H^1(A_K)$ are direct summands of $H^1 (KS(V))$
and $V\subset V_{1/2}\otimes H^1(A_K)$ (see \cite{geemen01} where the
other summands are also determined). The half twist is thus a partial
(it exists only when a CM field acts on $V$) generalization (it works
for any weight) of the Kuga-Satake construction.

The half twist is also implicit in the work of Kondo \cite{kondo99} on
$K3$ surfaces which are degree $4$ covers of $\bP^2$ totally branched
along plane quartics. These come with an action of $\bZ /4\bZ$ and the
transcendental part $V$ of $H^2 (S)$ is a vector space over the field
$\QQ(i)$. Kondo shows that the moduli space for such $V$ is isomorphic
to the moduli space of certain weight one Hodge structures $V'$ on
which the field $\QQ(i)$ acts.  In fact $V'$ is $V_{ 1/2}$, the half
twist of $V$. The moduli space for $V'$ is a quotient of a 6-ball and
Kondo proves that the moduli space of curves of genus 3 is birationally
isomorphic to this ball quotient. In the same paper, Kondo makes a
similar construction for the moduli space of curves of genus 4. Kondo's
work was motivated by results of Allcock-Carlson-Toledo
\cite{allcockcarlsontoledo98}, \cite{allcockcarlsontoledo00} who
produce a complex hyperbolic structure on the moduli space of cubic
surfaces by using cubic threefolds which are triple covers of $\bP^3$
totally branched along a cubic surface.

The half twist also appears in the study of variations of Hodge
structures. The families of weight two Hodge structures considered by
Carlson and Simpson in \cite{carlsonsimpson87} are obtained as
(negative) half twists of families of weight one Hodge structures with
a $\QQ(\sqrt{-1})$-action.

In this paper, we consider hypersurfaces $Y_k$ of degree $d\geq 3$ and
dimension $k\geq 1$ in $\bP^{ k+1 }$ which are $d$-fold covers of $\bP^k$
totally branched along hypersurfaces $X_{ k-1} $ of degree $d$ in
$\bP^k$. Such hypersurfaces come with an action of $\bZ/ d\bZ$. Letting
$\sigma$ be a generator of $\bZ/ d\bZ$, we define $V$ as the part of the
cohomology of $Y_k$ where the eigenvalues of $\sigma$ are primitive $d$-th
roots of $1$. Then $V$ is a vector space over the cyclotomic field $K$ of
$d$-th roots of unity. We will fix a CM-type $\Sigma_0$ of $K$ in
\ref{cycf} which is optimal for the examples under consideration. That is,
if $V$ has a half twist for some CM-type, then it has a half twist for this
CM-type.

We determine which $V$ allow a half twist.
If $V$ does not have
maximal level (i.e.\ if $V^{k,0}=0$), $V$ has a half twist. It
is then also interesting to see if $V(q)$, the Tate twist of $V$ whose
level is equal to its weight, has a half twist. In theorem \ref{lemht}
we determine the values of $d$ and $k$ for which the half twist of
$V(q)$ exists.

\

In case $V_{1/2}$ exists, we would like to have it as a
Hodge substructure of the cohomology of a `nice' algebraic variety. A
basic result on half twists (cf.\ Proposition \ref{halftens}) implies
that
\[
V_{1/2} (-1)\subset H^k(Y_k)\otimes
H^1(A_K)\subset H^{k+1}(Y_k\times A_K)
\]
where $A_K$ is an abelian variety with CM by $K$.
We find two nicer varieties:

\

\noindent
{\bf Theorem (cf. \ref{ZV}).} Let $Y_1$ be the Fermat curve of degree $d$.
\begin{enumerate}
\item There is an embedding $H^1 (A_K )\inj H^1 (Y_1)$. So if
$V_{1/2}$ exists, then
\[
V_{1/2} (-1)\subset H^k(Y_k)\otimes
H^1(Y_1)\subset H^{k+1}(Y_k\times Y_1).
\]
\item Let $Z_{ k+1 }$ be the $d$-fold cover of $\bP^{k+1 }$ totally
ramified along $Y_k$. If $V_{1/2}$ is exists, then 
\[
V_{1/2}(-1) \subset H^{k+1}_0(Z_{k+1}).
\]
\end{enumerate}

\
The second part of the above theorem follows from the first part and work
of Shioda.

Since $V_{1/2}(-1)$ is not of maximal level in $H^{k+1}_0(Z_{k+1})$ and
$H^{k+1}(Y_k\times Y_1)$, the general Hodge conjecture predicts that
$V_{1/2}$ is a Hodge substructure in the cohomology of subvarieties of
codimension one of $Z_{k+1}$ and $Y_k\times Y_1$. We don't know how to find
such a subvariety except in a few examples (cf.\ \cite{geemen01} for the
case of Calabi-Yau varieties).

\

A different and equally interesting type of question is a Torelli type
problem: to what extent does the half twist $V_{ 1/2}$ determine $X_
{ k-1 }$? In Section \ref{cub4} we prove

\

\noindent
{\bf Theorem (cf. \ref{genitor}).}
Suppose $d=3$, $k > 3$ and $V(q)_{1/2}$ is well-defined, i.e., by Theorem
\ref{lemht}, $k=3q+1$. Then the differential of the period map which
to $X_{ k-1 }$ associates $V(q)_{1/2}$ is generically injective.

\

An important application of our results is on Kuga-Satake
correspondances. For a Hodge structure $V$ of weight 2 with $V^{2,0}$
of dimension 1 there is the inclusion of Hodge structures $V\subset
H^1(KS(V))^{\otimes 2}$. Supposing $V$ is a Hodge substructure of the
cohomology of a variety $S$, the Hodge conjecture asserts that this
inclusion is induced by a correspondance on the product of $S$ and
$KS(V)^2$. We have

\

\noindent
{\bf Theorem (cf. \ref{KSthm}, \ref{KScor}).}
There exists a Kuga-Satake correspondance for the weight two Hodge
structure $V=H_0^4(Y_4)(1)$ where $Y_4$ is a triple cover of $\PP^4$
branched along a (general) smooth cubic threefold.

\

This generalizes a result of Voisin \cite{voisin96}. Theorem
\ref{KSthm} and Corollary \ref{KScor} below also apply to the quartic
surfaces considered by Kondo and may have applications to more general
situations involving the half twist.

\section{Polarized rational Hodge structures with automorphisms.}
\label{recalls}

\subsection{} We recall the basic results on half twists from \cite{geemen01}.
\subsection{Definition.} A (rational) Hodge structure of weight 
$k\;(\in\ZZ_{\geq 0})$
is a $\QQ$-vector space $V$ with a decomposition of its complexification 
$V_\CC:=V\otimes_\QQ\CC$ 
(where complex conjugation is given by 
$\overline{v\otimes z}:=v\otimes\bar{z}$ for $v\in V$ and $z\in\CC$):
$$
V_\CC=\oplus_{p+q =k}V^{p,q},\qquad\hbox{such that}\quad
\overline{V^{p,q}}=V^{q,p}, \qquad(p,\,q\in\ZZ_{\geq 0}).
$$
Note that we insist on $p$ and $q$ being non-negative integers throughout 
this paper, so we only consider `effective' Hodge structures.

\subsection{CM-type.}
Recall that a CM-field $K$ has $2r=[K:\QQ]$ complex embeddings
$K\hookrightarrow \CC$ which are pairwise conjugate and that a CM-type is a
subset $\{\sigma_1,\ldots,\sigma_r\}$ of distinct embeddings with the
property that no two are complex conjugate. Hence if we define embeddings
$\sigma_{r+i}(x):=\overline{\sigma_i(x)}$ then any embedding of $K$ in
$\CC$ is a $\sigma_j$ for some $j$, $1\leq j\leq 2r$.

\subsection{Half twists.}\label{hatw}
Let $V$ be a Hodge structure of weight $k$ on which a CM-field $K$ acts and
let $\Sigma=\{\sigma_1,\ldots,\sigma_r\}$ be a CM-type. The eigenspaces of
the $K$-action on the $V^{p,q}$'s are denoted by:
\[
V^{p,q}_j:=\{v\in V^{p,q}:\;xv=\sigma_j(x)v\quad \forall x\in K\},
\qquad 1\leq j\leq 2r.
\]
We define two subspaces of $V_\CC$ whose direct sum is $V^{p,q}$:
\[
V^{p,q}_+:=\oplus_{i=1}^r\,V^{p,q}_i,\qquad 
V^{p,q}_-:=\oplus_{i=1}^r\,V^{p,q}_{r+i}.
\]
We define the Hodge decomposition of the negative half twist of $V$
(w.r.t.\ $\Sigma$) by:
\[
V_{-1/2}^{p,q}:=V^{p-1,q}_+\oplus V^{p,q-1}_-.
\]
It is not hard to see that this is a Hodge structure of weight $k+1$
on $V$. By successively performing the negative half twist one obtains
$V_{-n/2}$, a Hodge structure on $V$ of weight $k+n$. However, half
twists do not give Tate twists: $V_{-2m/2}\neq V(-m)$. Half twists and
Tate twists are compatible in the sense that $(V_a)(b)=(V(b))_a$ for
all $a\in (1/2)\ZZ$ and $b\in\ZZ$.

A classical example is obtained by taking $V=K$, with the trivial
Hodge structure $V^{0,0}=V_\CC$. Then $K_{-1/2}\cong H^1(A_K,\QQ)$
for any abelian variety $A_K$ of dimension $r$ with CM by $K$ and
CM-type $\Sigma$.

\subsection{Positive half twists.} To define the (positive) half
twist one would put:
\[
V_{1/2}^{p,q}:=V^{p+1,q}_+\oplus V^{p,q+1}_-.
\]
This works if $V^{ k,0 } =0$. However, if $V^{ k,0 }\neq 0$, then the
subspaces $V^{k,0}_-$ and $V^{0,k}_+$ of $V_\CC$ do not appear in
$(V_{1/2})_\CC$ and therefore this definition does not define a Hodge
structure on $V$ (and in general not on any $\QQ$-subspace of $V$).
One can define the half twist of $V$ only if $V^{k,0}_-=0$ (the
complex conjugate of this space is $V^{0,k}_+$ which is then also
$0$). So one needs the eigenvalues of any $x\in K$ on $V^{k,0}$ to be
in the set $\{\sigma(x)\}_{\sigma\in\Sigma}$. The following
proposition shows that half twists appear naturally in certain tensor
products.

\subsection{Proposition.}(See \cite{geemen01}.)\label{halftens}
Let $V$ be a Hodge structure with CM by $K$ and fix a CM-type $\Sigma$ of
$K$.  Then we have an inclusion of Hodge structures (both of weight $k+1$):
\[
V_{-1/2}\subset V\otimes_\QQ K_{-1/2}, 
\]
given by:
\[
V_{-1/2}=\left\{
w\in V\otimes_\QQ K_{-1/2}:\; (x\otimes 1)w=(1\otimes x)w 
\quad\forall x\in K\right\}.
\]
Similarly, if $V$ admits a positive half twist, the Hodge structure
$V_{1/2}(-1)$ of weight $k-1+2=k+1$ is a Hodge substructure of
$V\otimes K_{-1/2}$:
\[
V_{1/2}(-1)=\left\{ w\in V\otimes_\QQ
K_{-1/2}:\; (x\otimes 1)w=(1\otimes \bar{x})w\quad\forall x\in K
\right\}.
\]

\section{Generalities on Hodge structures of hypersurfaces}\label{genHS}

\subsection{The Hodge substructure $V$.}\label{defV}\label{cycf}
For a smooth $(k-1)$-dimensional hypersurface $X_{ k-1 }$ of degree $d$
\[
X_{k-1}:=Zeroes(F_d(x_0,\ldots,x_k)) \qquad(\subset\PP^{k}),
\]
the cyclic $d$-fold cover of $\PP^k$ branched along $X_{k-1}$
is the smooth $k$-fold $Y_k$ defined by
\[
Y_k:=Zeroes(x_{k+1}^d+F_d(x_0,\ldots,x_{k}))\qquad(\subset\PP^{k+1}).
\]
We denote by $\alpha_k$ the automomorphism of order $d$ defined by
\[
\begin{array}{crclc}
\alpha_k: & Y_k &\longrightarrow & Y_k &\\
 & (x_0:\ldots:x_k:x_{k+1}) &\longmapsto & (x_0:\ldots:x_k:\zeta x_{k+1}) &
(\zeta=e^{ 2\pi i/d }).
\end{array}
\]
The action of $\alpha_k$ on $H^k_0(Y_k,\QQ)$ makes this
space a $\QQ[T]/(T^d-1)$-module. We will
denote by
\[ 
V\hookrightarrow H_0^k(Y_k,\QQ)
\]
the largest subspace on which the eigenvalues of $\alpha_k^*$ are
primitive $d$-roots of unity. If $d$ is a prime number we have
$V=H^k_0(Y_k,\QQ)$. The restriction of $\alpha_k^*$ to $V$ is denoted by
$\alpha$:
\[
\alpha:=(\alpha_k^*)_{|V}:V\longrightarrow V.
\]
The subring $K:=\QQ[\alpha]$ of $End(V)$ is isomorphic to the field of $d$-th roots of unity. The complex embeddings of $K$ are:
\[
\sigma_a:K\longrightarrow \CC,\qquad \sum x_j\alpha^j\longmapsto 
\sum x_j \zeta^{aj}\qquad(x_j\in\QQ ,\zeta = e^{ 2\pi i/d})
\]
with $a\in(\ZZ/d\ZZ)^*$. We define the CM-type $\Sigma_0$ for $K$ by:
\[
\Sigma_0:=\{\sigma_a\}_{0<a<d/2}.
\]
All half twists throughout this paper are taken w.r.t.\ this CM-type.
Thus $V$ has a half twist (w.r.t.\ $\Sigma_0$) if and only if all
eigenvalues of $\alpha$ on $V^{k,0}$ are in the set
$\{\sigma_a(\zeta)\}_{0<a<d/2}$.

We will often identify $K$ with $\QQ(\zeta)$ via $\sigma_1$. With this
convention, Proposition \ref{halftens} implies
that the half twist, when it exists, may also be defined via:
\[
V_{1/2}(-1)=\left\{ w\in V\otimes_\QQ
K_{-1/2}:\; (\alpha\otimes \zeta)w=w
\right\}
\]
since $1\otimes\bar{\zeta}=(1\otimes \zeta)^{-1}$.
We observe that if the half twist exists, then the eigenvalues of $\alpha\otimes\zeta$
on $V^{k,0}\otimes (K_{-1/2})^{1,0} $ are contained in
$\{\sigma_a(\zeta)\cdot\sigma_b(\zeta)\}_{0<a,b<d/2}$, so none of these is
equal to $1$ and hence $V_{1/2}$ is indeed an `effective' Hodge
structure of weight $k-1$.

We now review the basic results on Hodge structures of hypersurfaces
and the action of automorphisms on their cohomology and apply them to
$V$.

\subsection{Hodge numbers.}\label{secthn}
The Hodge numbers of a degree $d$ hypersurface $X_k$ in $\PP^{k+1}$ 
can be calculated using Griffiths residues
(see \cite{green94} page 44):
\[
H^{k-q , q} (X_k)\cong H^0\left(\frac{\cO_{\bP^{ k+1 }}}{ J} ( d (q+1)
  - k -2 )\right)
\]
where $J$ is the jacobian ideal of $X_k$, i.e., the ideal generated by 
the partial derivatives of an equation for $X_k$.

\subsection{Automorphisms.}\label{autos}
We compute the dimensions of the eigenspaces in the primitive
cohomology
\[
H^{k-q,q}_0(Y_k)(i):=\{x\in H^{k-q,q}_0(Y_k):\;\alpha_k^*x=\zeta^i x\;\}
\]
where $\alpha_k$ is the automorphism of $Y_k$ defined in \ref{defV} above.
The following lemma is an immediate consequence of a result of Shioda 
(cf.\ \cite{shioda79} Theorem I).

\subsection{Lemma.}\label{lemhpqi}
Let $J'$ be the Jacobian ideal of $F_d\in \CC[x_0,\ldots,x_{k}]$.
Then we have:
\[
h^{k-q,q}_0(Y_k)(i)= h^0\left(\frac{\cO_{\bP^{ k}}}{ J'} (d(q+1)
  - k-1-i)\right).
\]

\ts
According to \cite{shioda79}, Theorem I, the dimension of the
eigenspace is given by the number of $(k+1)$-tuples $a_0,\ldots,a_k$
with $1\leq a_j\leq d-1$ such that $(a_0+\ldots+a_k+i)/d=q+1$.
Moreover $h^{k-q,q}_0(Y_k)(i)=0$ if $i=0$ hence we take $1\leq
i\leq d-1$. Writing $b_j:=a_j-1$, so $0\leq b_j\leq d-2$, we get
$b_0+\ldots+b_k=d(q+1)-k-1-i$. To each such $(k+1)$-tuple we associate
the monomial $x_0^{b_0}\ldots x_k^{b_k}$, which is homogeneous of
degree $d(q+1)-k-1-i$ and has no variable to power $\geq d-1$. Since
these dimensions are constant in families of smooth hypersurfaces, we
may assume that $F_d=\sum x_i^d$ and we obtain the desired equality.
\qed

\subsection{The half twist for $V$.}\label{htv}
Since $V$ is the subspace of $H_0^k(Y_k,\QQ)$ on which the eigenvalues
of $\alpha_k^*$ are primititive $d$-th roots of unity, we have:
\[
V^{k-q,q}=\oplus_{i\in(\ZZ/d\ZZ)^*} \, H_0^{k-q,q}(Y_k)(i).
\]
In case $V^{k,0}=0$, the half twist of $V$ always exists. It may happen
that $V_{1/2}$ has higher level than $V$ (for example, a cubic surface
$Y_2$ has $V_\CC=V^{1,1}$ and thus level $0$, but $V_{1/2}$ has weight one
and hence has level $1$). In case $V^{k,0}=0$, we will therefore also
consider the Tate twist $V(q)$ of $V$ with $q$ choosen such that
$V^{k-q,q}\neq 0$ but $V^{u,s}=0$ if $u>k-q$. Then $V_{1/2}$ has lower
level than $V$ exactly when $V(q)$ has a half twist. Recall that all our
half twists are for the CM-type $\Sigma_0$ fixed in \ref{cycf}.

\subsection{Theorem.}\label{lemht}
Let $Y_k$ be a hypersurface of degree $d$ in $\PP^{k+1}$
as in Section \ref{defV}.

Define $q\in \ZZ_{\geq 0}$ and $t$ by:
\[
k=qd+t\qquad{\rm with}\quad t\in\,\{-1,\, 0,\,\ldots\,d-2\}.
\]
Then $H^{k-q,q}_0(Y_k)$ is the `extremal' summand of $H^k_0(Y_k,\CC)$, i.e.,
\[
H^{k-q,q}_0(Y_k)\neq 0,\qquad H_0^{u,s} (Y_k) =0\quad{\rm if} \;u>k-q.
\]
The Tate twist $V(q)$ of the Hodge substructure $V$ of $H_0^k(Y_k,\QQ)$ 
has a half twist if and only if
\begin{enumerate}
\item
either $d\not\equiv 2\;{\rm mod} \,4$ and:
\[
t >\frac{d-4}{2},
\]
\item 
or $d\equiv 2\;{\rm mod} \,4$ and
\[
t>\frac{d-6}{2}.
\]
\end{enumerate}

\ts  The integer $q$ is the smallest nonnegative integer for which $H^{ k-q
,q}_0 (Y_k)\neq 0$. Therefore, by the formula for the Hodge numbers of
$Y_k$ (Section \ref{secthn}), we have that $q$ is the nonnegative
integer such that $d( q+1 ) -k-2\geq 0$ and $dq -k -2 < 0$, i.e., $q$
is the integer satisfying
\[
q<\frac{k+2}{d}\leq q+1.
\]
Now it is clear that $k=qd+t$ as in the statement of the proposition.

We show first that $\dim H^{k-q,q}_0(Y_k)(i)\geq \dim H^{k-q,q}_0(Y_k)(i+1)$.
This implies that if $V$ has a half twist, then it has a half twist with
the CM-type $\Sigma_0=\{\sigma_a\}_{ 0<a<d/2 , a\in(\ZZ/d\ZZ)^* }$.
Put $a = d(q+1) -k-2=d-t-2$. Recall that
(Lemma \ref{lemhpqi})
\[
h^{k-q,q}_0(Y_k)(i)=
h^0\left(\frac{\cO_{\bP^{ k}}}{ J'} (a+1-i)\right).
\]
We may assume that $J'$ is generated by $x_0^{d-1},\ldots,x_k^{d-1}$, and
then $H^0((\cO_{\bP^k}/{ J'}) (a+1-i))$ is spanned by monomials
$x_0^{a_0}\ldots x_k^{a_k}$ with $\sum a_j=a+1-i$ and $a_j\leq d-2$ for all
$j$. Since $a=d-t-2$ and $-1\leq t\leq d-2$, multiplication by $x_k$
induces an injection
\[
H^{k-q,q}_0(Y_k)(i+1)\hookrightarrow H^{k-q,q}_0(Y_k)(i)
\]
except if $a=d-1$ and $i=1$. In this last case we can define a different
map which is also injective: we send $x_k^{d-2}$ to $x_0^{d-2}x_1$ and
we send all other monomials to their product by $x_k$. This inclusion
is not surjective except when $d=3$ and $k=2$, in which case $q=1$
and $V(-1)$, of weight $0$, does not have a half twist.

Next $V(q)$ has a half twist, for our choosen CM type $\Sigma_0$, if and
only if $H_0^{k-q,q}(Y_k)(i)=0$ for $i\in (\ZZ/d\ZZ)^*$ and $d/2<i<d$.
Using that these eigenspaces get smaller as $i$ gets bigger, we find that
$V(q)$ has a half-twist if and only if $h_0^{ k-q,q } (Y_k) (j) =0$ where
$j\in (\ZZ/d\ZZ)^*$ is the smallest unit with $j>d/2$.

If $d=2e$ with $e$ even or $d=2e+1$, the smallest unit in $(\ZZ/d\ZZ)^*$
larger than $d/2$ is $e+ 1\in (\ZZ/d\ZZ)^*$. Therefore $V$ has a half twist
if and only if $h_0^ { k-q,q } (Y_k) ( e+1 )=0$. This is equivalent to $a
-e <0$, that is, if $d=2e$ with $e$ even:
\[
(d-t-2)-d/2<0\;\Longleftrightarrow\; t>\frac{d-4}{2}
\]
and the same result holds if $d=2e+1$.

If $d=2e$ with $e$ odd, the smallest unit in $(\ZZ/d\ZZ)^*$ larger than
$d/2$ is $e+2$, hence $V$ has a half twist if and only if $a -e -1 <0$,
that is:
\[
(d-t-2)-(d/2)-1<0\;\Longleftrightarrow\;
t>\frac{d-6}{2}.
\]
\qed

\subsection{Corollary.}\label{cor}
Let $Y_k$ have degree $d$. Then $V\subset H^k(Y_k,\QQ)_0$ has a half
twist if and only if
\[
\begin{array}{crcl}
{\rm either}\quad & d < 2k+4  \quad&\quad {\rm and}&\quad k\equiv 0\;{\rm mod}\;2,\\
{\rm or}\quad & d \leq 2k+4\quad&\quad {\rm and}&\quad k\equiv 1\;{\rm mod}\;2.
\end{array}
\]

\ts If $d < k+2$, then $V^{ k,0 } = 0$ and, as we saw in \ref{htv}, $V$ has
a half twist. Else $k\leq d-2$ and we have $q=0$ and $t=k$. Therefore $V$
has a half twist if and only if either $k>(d-4)/2$, that is $d<2k+4$, or
$d\equiv 2\;{\rm mod} \,4$ in which case $k>(d-6)/2$ hence $d<2k+6$.  In
the last case we must actually have $d\leq 2k+4$ with equality only if $k$
is odd.  \qed

\section{The half twist of a cyclic cover}\label{Zh}

\subsection{Cyclic covers}\label{cyc}
In Section \ref{defV} we defined a Hodge substructure $V\subset
H^k_0(Y_k,\QQ)$ which has CM by the field $K$ of $d$-th roots of unity.
In case $V$ has a positive half twist $V_{1/2}$, we saw in Proposition
\ref{halftens} that we have an embedding
\[
V_{ 1/2 } (-1)\subset V\otimes K_{ -1/2 }\subset H^k_0(Y_k,\QQ)\otimes K_{
-1/2 }.
\]
Recall that $\alpha_k$ is the automorphism of
$Y_k=Zeroes(x_{k+1}^d+F_d(x_0,\ldots,x_k))$ defined by
\[
\alpha_k:Y_k\longrightarrow Y_k,\qquad (x_0:\ldots:x_k:x_{k+1})
\longmapsto (x_0:\ldots:x_k:\zeta x_{k+1})
\]
where $\zeta = e^{ 2\pi i/d}$. In particular, the Fermat curve $Y_1$
has the automorphism $\alpha_1$ (taking $F_d:=x_0^d+x_1^d$).\\ The
following lemma shows that we have a geometric realization of the half
twist:
\[
V_{ 1/2 } (-1)\subset H^k_0 (Y_k ,\QQ )\otimes H^1 (Y_1 ,\QQ )\subset H^{
k+1 } (Y_k\times Y_1 ).
\]

\subsection{Lemma.}\label{lemZV}
The Hodge structure $K_{-1/2}$ is a Hodge substructure of
$H^1(Y_1,\QQ)$.

\ts
On the Fermat curve $Y_1$ defined by $x_2^d+x_1^d+x_0^d$ consider 
the automorphism 
\[
\gamma:Y_1\longrightarrow Y_1,\qquad 
(x_0:x_1:x_2)\longmapsto (x_0:\zeta^{-1}x_1:\zeta x_2).
\]
The holomorphic one forms on $Y_1$ are the $y_1^{a-(d-1)}y_2^b{\rm
d}y_2$ (where we put $y_i =\frac{ x_i }{ x_0 }$), with $0\leq a+b\leq
d-3$ and $a,\,b\geq 0$. The $\gamma$-invariants have
$-a+(d-1)+b+1\equiv 0$ mod $d$, that is $a\equiv b$ mod $d$. Hence the
possible values for $(a,b)$ are $(0,0),\,(1,1),\ldots,([\frac{ d-3 }{
2}],[\frac{ d-3 }{ 2}])$. In particular $h^1(Y_1,\QQ)^\gamma=
2[\frac{ d-1 }{ 2}]$ and the eigenvalues of $\alpha_1^*$ on 
$H^{ 1,0 }( Y_1)^{\gamma}$ are
$\zeta,\,\zeta^2,\,\ldots,\zeta^{[ (d-1) /2] }$. Since $\alpha_1^d=1$
and $\alpha_1$ commutes with $\gamma$,
the space $H^1(Y_1,\QQ)^\gamma$ is a $\QQ[T]/(T^d-1)$-module (via
$T\mapsto\alpha_1$). The eigenvalue $\zeta$ has multiplicity one,
hence $H^1(Y_1,\QQ)^\gamma$ has a unique submodule $M$ isomorphic to
$K = \QQ (\zeta)$. As a Hodge structure, $M\cong K_{-1/2}$ since
both have the same CM-type. This concludes the proof of the lemma.
\qed

\subsection{}
In fact we have a better result which will allow us (see below) to embed
$V_{1/2} (-1)$ in $H^{k+1}(Z_{k+1},\QQ)$ where
\[
Z_{k+1}:=Zeroes(x_{k+2}^d+x_{k+1}^d+F_d(x_0,\ldots,x_{k}))
\qquad(\subset\PP^{k+2}).
\]

\subsection{Theorem.}\label{ZV}
With the notation of \ref{defV} and \ref{cyc}, assume that the half twist
$V_{1/2}$ of the Hodge structure $V\;(\subset H^k_0(Y_k,\QQ)$) exists. Let
$W$ be the Hodge substructure of $H^k_0(Y_k,\QQ)\otimes H^1(Y_1,\QQ)$ on
which the automorphism $\beta:=(\alpha_k^*,\alpha_1^*)$ acts trivially:
\[
W :=\left(H^k_0(Y_k,\QQ)\otimes H^1(Y_1,\QQ)\right)^{\langle\beta\rangle}.
\]
Then there is an inclusion of Hodge structures of weight
$k+1$:
\[
V_{1/2}(-1)\hookrightarrow W =
\left(H^k_0(Y_k,\QQ)\otimes H^1(Y_1,\QQ)\right)^{\langle\beta\rangle}.
\]

\ts We have that 
\[
W = \left( V\otimes K_{ -1/2 }\right)^{\langle\beta\rangle }\oplus\ldots .
\]
Since $\beta(v)=(\alpha_k^*\otimes\alpha_1^*)v$ for $v\in
H^k_0(Y_k,\QQ)\otimes H^1(Y_1,\QQ)$ we also have:
\[
\beta(v)=v \qquad\hbox{if and only if}\qquad (\alpha_k^*\otimes 1)v=(
1\otimes\bar{\alpha_1^* })v
\]
hence, by Proposition \ref{halftens}, we have $\left( V\otimes
K_{-1/2}\right)^{\langle\beta\rangle}\cong V_{1/2}(-1)$
and thus $V_{1/2}(-1)\subset W$. \qed

\

We have the following immediate consequence of results of Shioda.

\subsection{Proposition.}(Shioda.)
\label{cohY}
The primitive cohomology of $Z_{k+1}$ has a direct sum decomposition
\[
H^{k+1}_0(Z_{k+1},\QQ)\cong H^{k-1}_0(X_{k-1},\QQ)(-1)^{\oplus
  (d-1)}\,\oplus W.
\]

\ts In the case of Fermat varieties the direct sum decomposition of $H_0^{
k+1 } (Z_{ k+1 } ,\QQ)$ is proved in \cite{shioda79}, Theorem II and it
applies in this case as well since we deal with a topological property.
\qed

\subsection{Corollary}\label{corV}
We have an inclusion of Hodge structures of weight $k-1$:
\[
H_0^{k-1}(X_{k-1},\QQ)^{\oplus (d-1) }\oplus V_{1/2}
\;\hookrightarrow\;
H_0^{k+1}(Z_{k+1},\QQ)(1).
\]

\ts
Immediate consequence of \ref{cohY} and \ref{ZV}.
\qed

\

\noindent
We could then ask whether we can have
\[
V_{ 1/2 } (-1) = W.
\]
The following dimension computation shows that this is the case if and only
if $d=3$.
\subsection{Lemma.}\label{lemdim}
For all $k\geq 2$ we have the identity
\[
 h^{ k+1 }_0 = (d-1) h^{k-1}_0 + (d-2) h^k_0.
\]
where $h^k_0$ is the dimension of the primitive cohomology of a smooth
$k$-dimensional hypersurface of degree $d$ in $\PP^{ k+1 }$. Therefore, by
Proposition \ref{cohY}
\[
\dim W=(d-2) h^k_0.
\]

\ts It suffices to prove that $ h^{ k+1 }_0 = (d-1) h^{k-1}_0 + (d-2)
h^k_0$ for special hypersurfaces $X_k$ which are cyclic degree $d$ covers
of $\bP^k$ totally ramified along a degree $d$ hypersurface $X_{
k-1}$. Then the Zeuthen-Hurwitz formula gives
\[
\chi_{top} (X_k )= d\chi_{ top} (\bP^k ) - (d-1)\chi_{top } (X_{ k-1 })
\]
or
\[
k+1 + (-1)^k h^k_0 = d(k+1 ) -(d-1) (k + (-1)^{ k-1 } h^{ k-1}_0)
\]
or
\[
(-1)^k h^k_0 = (d-1) ( 1- (-1)^{k-1} h^{ k-1 }_0).
\]
From this we calculate
\[
\begin{array}{rcl}
(-1)^{ k+1 } h_0^{ k+1 } + (-1)^k (d-2) h^k _0 &=& (d-1) ( 1 - (-1)^k
h^k_0 ) + (-1)^k (d-2) h^k_0\\
&= & (d-1) - (-1)^k h^k_0\\
& =& (d-1) - (d-1) ( 1 - (-1 )^{ k-1} h^{ k-1}_0 )\\
& =& (d-1) (-1)^{ k-1} h^{ k-1 }_0,
\end{array}
\]
multiplying by $(-1)^{ k+1}$ gives the desired result.
\qed

\subsection{Cohomology of cubics.}\label{cubtw}
In case $Y_k$ has degree $3$, we have $V=H^k_0 (Y_k,\QQ)$. Since for $k>1$
we have $V^{k,0}=0$, its half twist exists and \ref{cohY}, \ref{ZV} and
\ref{lemdim} imply that
\[
H^{k+1}_0(Z_{k+1},\QQ)(1)=
H^{k-1}_0(X_{k-1},\QQ)^{\oplus 2}\,\oplus\, H^k_0(Y_k,\QQ)_{1/2}.
\]
To see for which values of $k$ we can apply a half twist to
$V(q)$, we apply Theorem \ref{lemht} and find that $k=
3q+1$. In that case we have $h^{k-q,q}_0(Y_k)=1$.

\subsection{Cohomology of Quartics.} 
In case $Y_k$ has degree 4, Theorem \ref{lemht} with $d=4$ shows that the
Hodge substructure $V(q)\subset H^k(Y_k,\QQ)$ has a half twist if and only
if $k=4q+1$ or $4q+2$. In case $k=4q+2$ we have $\dim
V^{k-q,q}=h^{k-q,q}_0(Y_k)=1$. From Corollary \ref{cor} it follows that
$V$ has a half twist for any $k$ (in fact $V^{k,0}=0$ for $k>2$).
We have:

\subsection{Lemma.}\label{lemquartics} 
The Hodge structure $W$ splits as follows:
\[
W(1) =
V_{1/2}^{\oplus 2}\oplus \left(V'\otimes K_{ -1/2 }\right),
\]
where $V'\subset H^k_0(Y_k,\QQ)$ is the Hodge substructure on which
$\alpha_k$ acts as $-1$. Therefore by \ref{cohY}
\[
H^{k+1}_0(Z_{k+1},\QQ)(1)\cong H^{k-1}_0(X_{k-1},\QQ)^{\oplus 3}\oplus
V_{1/2}^{\oplus 2}\oplus \left(V'\otimes K_{ -1/2 }\right).
\]

\ts We have
$H^k_0(Y_k,\QQ)=V\oplus V'$ and $H^1(Y_1,\QQ)\cong K_{-1/2}^{\oplus
  3}$ because the action of $\alpha_1^*$ on $H^{1,0}(Y_1)$ has
eigenvalue $i$ with multiplicity $2$ and eigenvalue $-1$ with
multiplicity $1$ (the $(-1)$-eigenspace is a copy of $K_{-1/2}$, using
permutations of the coordinates one finds the other two copies). Hence:
\[
W=\left(H^k_0(Y_k,\QQ)\otimes H^1(Y_1,\QQ)\right)^{\langle\beta\rangle}=
V_{1/2}^{\oplus 2}\oplus \left(V'\otimes K_{-1/2}\right).
\]
\qed

\section{Surfaces.}

\subsection{} We consider surfaces $Y_2\subset \PP^3$ of degree $d$.
The Hodge substructure $V\subset H^2(Y_2,\QQ)$ has a half twist if and only
if $d<8$ (see Corollary \ref{cor}). We briefly discuss various cases.

\subsection{Quartic surfaces.}\label{kondo}
We study the case $d=4$ in some detail since it is related to
recent work of Kondo \cite{kondo99}. Start with a curve
$C:=X_1\;(\subset \PP^2)$ of degree 4. The surface $Y_2$ is a K3
surface, in particular $H^{2,0}(Y_2)=1$. The primitive cohomology of
a general $Y_2$ decomposes into the direct sum of its transcendental
part $V$ and its algebraic part $NS_0 (Y_2 )_{\QQ }$:
\[
H^2_0(Y_2,\QQ)=V \oplus NS_0(Y_2)_\QQ,\qquad \dim_\QQ V=14,\qquad
\dim NS_0 (Y_2 )_{\QQ }=7.
\]
Kondo \cite{kondo99} proves, (using the Torelli theorem for
K3 surfaces) that the map $C\mapsto V_\ZZ$ defines an injective
morphism from the moduli space of non-hyperelliptic curves of genus 3
to a 6-ball quotient and studies its extension
to hyperelliptic curves of genus 3 and to nodal plane quartics.
This ball quotient is a moduli space of abelian varieties of
dimension 7 with (1,6)-action by $K=\QQ(i)$.

Kondo's construction can also be done with half-twists: We can apply
the half twist to the transcendental part $V$. The
Hodge structure $V_{1/2}$ is of weight one, has dimension 14 and
defines an (isogeny class of) abelian variety $A_C$ of dimension 7 on
which the action of the field $K$ is of type (1,6):
\[
V_{1/2}=H^1(A_C,\QQ),\qquad \dim A_C=7.
\]

Since we can also carry out the half twist over $\ZZ$, that is we can
define $(V_\ZZ)_{1/2}$ by the same formula as for $V_{1/2}$, we
can in fact define an abelian variety $A_C$ (with
$H^1(A_C,\ZZ)=(V_\ZZ)_{1/2}$ and polarization induced by the
polarization $\psi$ on $V_\ZZ$ via $E(v,w):=\psi(v,iw)$). We can of
course recover $V_\ZZ$ from $H^1(A_C,\ZZ)$. This abelian variety is
in fact Kondo's abelian variety.

Since $V_{1/2}\subset H^3(Z_3,\QQ)$, Kondo's abelian variety
is isogeneous to a subvariety of the intermediate jacobian of the quartic 
threefold $Z_3$.
It is well-known that
\[
H^3(Z_C,\CC)=H^{2,1}(Z_C)\oplus H^{1,2}(Z_C),\qquad {\rm and}\quad
h^{2,1}(Z_C) =30.
\]
In particular, the intermediate Jacobian $J_C$ of $Z_C$ is a 
(principally polarized) abelian variety.  Lemma \ref{lemquartics} implies that 
$$
J_C\sim J(C)^3\times A_C^2\times A_K^7
$$
since now $V'\cong \QQ(-1)^7\subset NS_0(Y_2)_\QQ$ is a trivial
Hodge structure. Here $A_K$ is an elliptic curve with CM by $K$.

\subsection{Quintic surfaces.}
\label{quintics}
Since $d$ is prime we have $V=H^k_0(Y_k,\QQ)$.  Theorem \ref{lemht}
implies that $V$ has a half twist if and only if $k=5q+2$ or $k=5q+3$.
In the last case $h^{k-q,q}_0=1$, but if $k=5q+2$ then:
\[
h^{k-q,q}_0(Y_k)=
h^0\left(\frac{\cO_{\bP^{ k+1 }}}{ J} (5(q+1)-(5q+2) -2)\right)
= h^0\left(\cO_{\bP^{ k+1 }} (1)\right)
=k+2
\]
and the eigenspaces have dimensions (cf.\ Lemma \ref{lemhpqi}):
\[
h_0^{k-q,q}(1)= h^0(\cO_{\bP^{ k}}(1))=k+1,\qquad
h_0^{k-q,q}(2)= h^0(\cO_{\bP^{ k}}(0))=1.
\]

We restrict our attention to a quintic surface $Y_2$, a 5:1 cover of
$\PP^2$ branched along a plane quintic. Then $V_{1/2} (-1)\subset
H^3(Z_3,\QQ)$, but since $H^{3,0}(Z_3)\neq 0$ we cannot expect to see
$V_{1/2}$ easily. However, the folowing geometric construction seems to
induce the inclusion $V=H^2_0(Y_2,\QQ)\subset V_{1/2}\otimes H^1(A_K)$.

Consider a general line $l$ in $\PP^2$. Its inverse image under the cover
$Y_2\rightarrow \PP^2$ is isomorphic to a plane curve $C=C_l$ defined by an
equation of type $y^5=f_5(x)$ from some quintic polynomial
$f_5\in\CC[x]$. The moduli space of such curves is 2-dimensional. We can
decompose the cohomology of $C$ with respect to the automorphism of order 5
(cf.\ \ref{lemhpqi}) and obtain:
\[
h^{1,0}(C)(1)=3,\qquad h^{1,0}(C)(2)=2,\qquad h^{1,0}(C)(3)=1,\qquad
h^{1,0}(C)(4)=0.
\]
The moduli space of abelian 6-folds with this type of automorphism, let's
call it ${\cal M}$, also has dimension 2 (since the space of invariants for
the automorphism in $S^2H^{1,0}(C)$ is also 2-dimensional, see also
\cite{jongnoot91} for this example). In ${\cal M}$ we have curves
parametrizing abelian varieties isogeneous to products $A_2\times A_4$
where $A_2$ is the jacobian of $C_2:\;y^2=x^5-1$ (and
$h^{1,0}(C_2)(1)=h^{1,0}(C_2)(2)=1$, the other $h^{1,0}(i)$'s being zero)
and $A_4$ is an abelian 4-fold with automorphism of order 5 with
eigenvalues $\zeta,\,\zeta,\,\zeta^2,\,\zeta^3$ on $H^{1,0}(A_4)$. In
particular, we can take $A_K=A_2$.

We have a rational map 
\[
(\PP^2)^*\longrightarrow {\cal M},\qquad
l\longmapsto Jac(C_l).
\]
The inverse image $Z$ of a curve in ${\cal M}$ as above carries a family of 
curves ${\cal C}$ which maps onto $Y_2$:
\[
\begin{array}{rcl}
{\cal C}&\longrightarrow &Y_2\\
\downarrow&&\\
Z
\end{array}\qquad\qquad
\begin{array}{rcl}
{\cal C}_l=C_l&\hookrightarrow &Y_2\\
\downarrow&&\\
\{ l\}
\end{array}
\]
The Jacobians of the $C_l$'s have a constant factor $A_2$ and one
expects that the pull-back of $V$ lies in $H^1(A_2,\QQ)\otimes
H^{1}(Z,\QQ)$ and that $V_{1/2}\subset H^1(Z,\QQ)$.

\subsection{Sextic surfaces.}
The surface $Y_2$ is a 6:1 cover of $\PP^2$ branched along a plane sextic
curve. The action of the automorphism $\alpha_2$ splits the 105-dimensional
primitive cohomology in 3 subspaces:
\[
H_0^2(Y_2,\QQ)=V_6\oplus V_3\oplus V_2,\qquad V:=V_6,
\]
where the eigenvalues of $\alpha_k^*$ on the $V_i$ are primitive $i$-th
roots of unity and, by definition, $V=V_6$. From Proposition
\ref{lemht} we see that $V$ has a half twist. From the formulae from
section \ref{genHS} it follows that
\[
\dim V=42,\quad
\dim V^{2,0}= h^{2,0}(Y_k)(1)=6,
\]
\[
\dim V^{1,1}= h^{1,1}(Y_k)(1)+ h^{1,1}(Y_k)(5)=15+15=30.
\]

Moreover, in this case $V_2$, which has CM by the field of cube
roots of unity, also has a half twist:
\[
\dim V_2=42,\quad
\dim V_2^{2,0}= h^{2,0}(Y_k)(2)=3,
\]
\[
\dim V_2^{1,1}= h^{1,1}(Y_k)(2)+ h^{1,1}(Y_k)(4)=18+18=36.
\]
We do not know a geometrical interpretation for the half twists of
$V_2$ and $V_3$.

\subsection{Septic surfaces.}
Since $7$ is a prime number, $V=H_0^2(Y_2,\QQ)$. The main problem here
is to find a correspondance which induces the inclusion
$V\hookrightarrow V_{1/2}\otimes K_{-1/2}$. In this case we can take
$K_{-1/2}=H^1(C,\QQ)$ where $C:\;y^2=x^7-1$ ($\zeta\in K$ acts on the
curve via $(x,y)\mapsto (x,\zeta y)$). Since $V_{1/2}$ has weight one
it may be identified with a Hodge substructure of the $H^1$ of some
curve $C'$. Of course, we also have $V_{1/2}(-1)\subset H^3(Z_3,\QQ)$
but it seems difficult to describe this Hodge substructure
geometrically (via an abel-jacobi map).

The Hodge conjecture predicts the existence of a cycle on $Y_2\times
(C\times C')$ which induces an inclusion $H^2_0(Y_2,\QQ)\hookrightarrow
H^1(C,\QQ)\otimes H^1(C',\QQ)$.

\section{Kuga-Satake correspondances}\label{KS}

\subsection{Kuga-Satake varieties: general.}
For a polarized Hodge structure $V$ of weight 2 with $\dim V^{2,0}=1$
Kuga and Satake \cite{kugasatake67} defined an abelian variety
$J_{KS}(V)$ with the property that
\[
V\hookrightarrow H^1(J_{KS}(V),\QQ)^{\otimes 2}.
\]
Voisin \cite{voisin96} observed that if $V$ is of CM-type for an
imaginary quadratic field $K$, then $J_{KS}$ has abelian subvarieties
of dimension $1$ and $n$, where $\dim_\QQ V=2n$. More precisely (cf.\ 
\cite{geemen01}), there is an inclusion of Hodge structures:
\[
K_{-1/2}\oplus V_{1/2}\hookrightarrow H^1(J_{KS},\QQ).
\]
For such a Hodge structure $V$ the inclusion of $V$ into
$H^1(J_{KS}(V),\QQ)^{\otimes 2}$ is now simply a consequence of
Propositon \ref{halftens} which gives an inclusion:
\[
V\hookrightarrow V_{1/2}\otimes K_{-1/2} \qquad(\subset
H^1(J_{KS}(V),\QQ)\otimes H^1(J_{KS}(V),\QQ)).
\]
In case $V\subset H^{2k}(Y,\QQ)(k-1)$, a basic problem in Hodge theory
is to find a cycle, a Kuga-Satake correspondance, on the product of
$Y$ and $J_{KS}(V)^2$ which induces this inclusion of Hodge
structures. The existence of such a cycle is predicted by the Hodge
conjecture. The following theorem shows that such a cycle exists for
example in the case where $V$ is the primitive cohomology of a cubic
fourfold $Y=Y_4$; the case of a quartic surface $Y_2$ was already
considered in \cite{geemen01}.

\subsection{Theorem.}\label{KSthm}
In the cases $d=3$, $k=4$ and $d=4$, $k=2$, there is a cycle $U$ on the
product of $Y_k$ and $Y_k\times Y_1\times A_K$ such that the map
\[
[U]:H^k(Y_k,\QQ)\longrightarrow 
H^{k+1}(Y_k\times Y_1,\QQ)\otimes H^1 (A_K ,\QQ)\quad(\subset
H^{k+2}(Y_k\times Y_1\times A_K,\QQ))
\]
induces
\[
V\stackrel{\cong}{\longrightarrow} V(-1)\;\subset 
V_{1/2}(-1)\otimes K_{-1/2}
 \subset  H^{k+1}(Y_k\times Y_1 ,\QQ)\otimes H^1 (A_K ,\QQ)
\]
where $V\subset H_0^k(Y_k,\QQ)$ is as in \ref{defV}.

In particular, $U$ is a Kuga-Satake correspondance.

\ts In these cases the Hodge substructure $V\subset H^k(Y_k,\QQ)$
has level two and, after a Tate twist in case $k=4$, we have $\dim
V^{2,0}=1$. Moreover $V$ is a vector space over $K$. 

Let $\gamma$ be the automorphism of the Fermat curve $Y_1$ which
appears in the proof of Lemma \ref{lemZV}. If $d=3,\,4$ then
$Y_1/\gamma\cong A_K$ where $A_K$ is an elliptic curve with CM by $K$,
i.e., $H^1 (A_K,\QQ )\cong K_{-1/2}$. So if $f:Y_1\rightarrow A_K$ is the
quotient map then $ f^* H^1 (A_K ,\QQ) = H^1 (Y_1 ,\QQ)^\gamma$. Let
$\mu_f$ be the composition of $1\otimes f^*$ and the projection onto the
space of $\beta$-invariants:
\[
\mu_f: H^k(Y_k ,\QQ)\otimes H^1(A_K ,\QQ)\longrightarrow (H^k(Y_k
,\QQ)\otimes H^1(Y_1 ,\QQ))^{\langle\beta\rangle} = W\subset H^k(Y_k
,\QQ)\otimes H^1(Y_1 ,\QQ).
\]
Then $V\otimes K_{ -1/2}$ is mapped onto $V_{1/2}(-1)$ and $\mu_f$
is induced by a correspondance.

Now consider the map:
\[
H^k(Y_k)\otimes K_{ -1/2}\otimes K_{ -1/2}\stackrel{\mu_f\otimes
  1}{\longrightarrow} W\otimes K_{ -1/2}\subset H^k(Y_k
,\QQ)\otimes H^1(Y_1 ,\QQ)\otimes K_{ -1/2 }.
\]
It maps $V\otimes K_{ -1/2} \otimes K_{ -1/2}$ onto $V_{1/2}(-1)\otimes
K_{ -1/2}$.
We define a Hodge substructure of $H^k(Y_k)\otimes K_{ -1/2}\otimes
K_{ -1/2}$ by:
\[
S:=\{w\in V\otimes K_{ -1/2}\otimes K_{ -1/2}:\;
(\alpha\otimes\zeta\otimes 1)w=w,\quad (1\otimes\zeta\otimes\zeta)w=w\,\}.
\]
Obviously we have $S\subset V_{1/2}(-1)\otimes K_{ -1/2}$ and
therefore $\mu_f$ embeds $S$ into $W\otimes K_{ -1/2}$. Since the
$(\zeta\otimes\zeta)$-invariant subspace of $K_{-1/2}\otimes K_{-1/2}$
is $K(-1)$, a trivial Hodge structure of weight 2 and rank 2, we have
that $S\subset V(-1)^{\oplus 2}$.  For dimension reasons it follows
that $S\cong V(-1)$ and thus $(\mu_f\otimes 1)(S)\cong V(-1)$. Since
the Hodge substructure $K(-1)\subset K_{ -1/2}\otimes K_{ -1/2}\cong
H^1(A_K)\otimes H^1(A_K)$ is trivial there is an algebraic cycle
$\Gamma\subset Y_k\times (Y_k\times A_K^2)$ such that
$[\Gamma]:H^k(Y_k)\longrightarrow H^{k+2}(Y_k\times A_K^2)$ induces
the (twisted) isomorphism $V\rightarrow S\cong V(-1)$.

The desired correspondance is the composition of $\Gamma$ and
$\mu_f\otimes 1$. On cohomology and restricted to $V\subset H^k(Y_k,\QQ)$
one has:
\[
V\stackrel{\Gamma}{\longrightarrow} V\otimes K_{ -1/2}\otimes K_{
  -1/2} \;\stackrel{\mu_f\otimes 1}{\longrightarrow}\;
V(-1)\qquad(\subset H^{k+1}(Y_k\times Y_1)\otimes H^1(A_K)).
\]
\qed

\subsection{Corollary.}\label{KScor}
In the cases $d=3$, $k=4$ and $d=4$, $k=2$, there is a cycle $U_Z$ on the
product of $Y_k$ and $Z_{k+1}\times A_K$ such that the map
\[
[U_Z]:H^k(Y_k,\QQ)\longrightarrow 
H^{k+1}(Z_{k+1},\QQ)\otimes H^1 (A_K ,\QQ)\quad(\subset
H^{k+2}(Z_{k+1}\times A_K,\QQ))
\]
induces
\[
V\stackrel{\cong}{\longrightarrow} V(-1)\;\subset 
V_{1/2}(-1)\otimes K_{-1/2}
 \subset  H^{k+1}(Z_{k+1},\QQ)\otimes H^1 (A_K ,\QQ)
\]
where $V\subset H_0^k(Y_k,\QQ)$ is as in \ref{defV}.

In particular, $U_Z$ is a Kuga-Satake correspondance.

\ts In the cases under consideration, $H_0^{ k+1 }(Z_{k+1},\QQ)$ is a Hodge
structure of level one (so, after a Tate twist, it is of weight one). Thus
the intermediate jacobian of $Z_{k+1}$ is an abelian variety and
$V_{1/2}(-1)$ defines an abelian subvariety.

The work of Shioda (\cite{shioda79} Theorem II) provides a correspondance
(a blow up of $Y_k\times Y_1$) which induces a map (as in Proposition
\ref{cohY})
\[
\mu: H^k(Y_k)\otimes H^1(Y_1)\longrightarrow H^{k+1}(Z_{k+1},\QQ).
\]
This map embeds the subspace $\left(H^k(Y_k)\otimes H^1(Y_1)\right)^\beta$
into $H^{k+1}(Z_{k+1},\QQ)$.

The composition of $[U]$ in the previous theorem and $\mu\otimes 1$ is
a correspondance $[U_Z]$ as required.
\qed

\subsection{Remark.}\label{voisin}
The cohomology of a cubic fourfold $Y$ with equation
$x_5^3+x_4^3=G(x_0,\ldots,x_3)$ was considered by Voisin
\cite{voisin96}, she shows that a Kuga-Satake correspondance exists
using the variety of lines on the cubic threefold defined by
$x_4^3=G(x_0,\ldots,x_3)$.  Our result is a generalization of hers.
See also \cite{geemen01} for Kuga-Satake varieties of weight two Hodge
structures on which an imaginary quadratic field acts.

\section{Cubic fourfolds and the geometry of cubic cyclic covers.} \label{cub4}

\subsection{}
We again consider a cubic 4-fold $Y_4\subset\PP^5$ with equation
$x_5^3=F_3(x_0,\ldots,x_4)$ where $F=0$ defines a smooth cubic 3-fold
in $\PP^4$. We already obtained an interesting result on the
cohomology of $Y_4$ in Theorem \ref{KSthm} and Corollary \ref{KScor}.

The Hodge numbers of $Y_4$ are:
\[
 h^{4,0}(Y_4)=0,\qquad h^{3,1}(Y_4)=1,\qquad h^{2,2}_0(Y_4)=20.
\]
The Hodge structure $V:=H^4_0(Y_4,\QQ)(1)$ has weight 2 and
$h^{2,0}=1$. We can define a half twist of $V$
\[
 V:=H^4_0(Y_4,\QQ)(1),\qquad V_{1/2}:=H^4_0(Y_4,\QQ)(1)_{1/2}. 
\]
Thus $V_{1/2}$ is a 22-dimensional polarized weight 1 Hodge structure
of CM-type with field $K$. Let $J(V_{1/2})$ be an abelian
variety such that
\[
V_{1/2}\cong H^1(J(V_{1/2}),\QQ).
\]
Note that $\dim J(V_{1/2})=11$ and that $K\subset End_\QQ(J(V_{1/2}))$.

The
middle cohomology of the cubic 5-fold $Z_5$ has level 1, 
that is, the intermediate
Jacobian of $Z_5$ is an abelian variety. By \ref{ZV}, \ref{cohY} and
\ref{lemdim} the 21-dimensional intermediate jacobian $J(Z_5)$
of the cubic five-fold $Z_5$ is isogeneous to 2 copies of the
5-dimensional intermediate Jacobian $J(X_3)$ of the cubic 3-fold $X_3$
and one copy of $J(V_{1/2})$:
\[
J(Z_5)\sim J(X_3)^2\times J(V_{1/2}).
\]

\subsection{Lemma.}
The eigenvalues of $x\in K$ on the tangent space to $J(V_{1/2})$ at
the origin are $x$ with multiplicity 1 and $\bar{x}$ with multiplicity
10.

\ts
By definition, we have $V_{1/2}^{1,0}=H^{3,1}_0(Y_4,\CC)_+\oplus
H^{2,2}_0(Y_4,\CC)_-$ on which $x\in K$ acts via diagonal matrices
$diag(\sigma(x),\bar{\sigma}(x))$. Since $V^{1,0}$ is the dual of the
tangent space to the origin, the lemma follows if we choose the right
embedding $K\subset\CC$.
\qed 

\subsection{Ball quotients.}
The moduli space of Hodge structures of weight one and of type (1,10) for
the field $K$ is the 10-ball. Also cubic threefolds depend on 10
moduli. By Lemma \ref{genitor} below, a general cubic threefold is sent
to a general abelian variety in this 10-dimensional family which is simple,
and, in particular, it is not (directly) related to the 5-dimensional
intermediate jacobian of the cubic 3-fold.

In this context one needs to work with $\ZZ$-Hodge structures.  One
can take $H^1(J(V_{1/2},\ZZ)) = H^4_0(Y_4,\ZZ)$, and the polarization
$E$ on $H^1(J(V_{1/2},\ZZ))$ is determined by the polarization $Q$ on
$H^4_0(Y_4,\ZZ)$, which was determined by Hassett \cite{hassett98}
following work of Beauville and Donagi \cite{beauvilledonagi85}, and
the action of the order three automorphism of $Y_4$ on the polarized
Hodge structure $H^4_0(Y_4,\ZZ)$. Unfortunately, this action (i.e.\ 
its conjugacy class in $O(Q)(\ZZ)$) is not known.

Independently, Allcock, Carlson and Toledo have extended their results
(private communication): they can show that the moduli spaces of cubic
threefolds and that of all Del-Pezzo surfaces are ball
quotients. For cubic threefolds, they show that the moduli space of
the Hodge structures of $H^4 (Y_4)$ is an arithmetic quotient of a
$10$-ball (if, instead, one uses half-twists as above, this is immediate
since $V_{1/2}$ has weight $1$) and (using Voisin's Torelli theorem
\cite{voisin86}) that the period map $X_3\mapsto H^4 (Y_4)$ is
injective. The case of Del-Pezzo surfaces of degree $2$ is in the
paper \cite{carlsontoledo99} by Carlson and Toledo and is equivalent to
Kondo's construction. The aim of Allcock, Carlson and Toledo is to use
the ball quotients to determine the fundamental groups of the
corresponding moduli spaces and their monodromy representations.

\subsection{Geometry of cubic cyclic covers.}
We consider again  a cubic $k$-fold $Y_k$.
As we observed before (see \ref{htv}), $V$ has a half twist for all
$k$, but its level will be one more than that of $V$ except if
$k=3q+1$ in which case it will be one less (cf. \ref{cubtw}).
In particular, for any $k$, there is an inclusion 
\[
V\hookrightarrow V_{ 1/2 }\otimes K_{ -1/2 }.
\]

It is easy to construct a variety $B$ of dimension $k-1$ and a
dominant rational map $B\times A_K\longrightarrow Y_k$ where $A_K$ is,
as before, an elliptic curve with CM by $K$, i.e., $H^1 (A_K)\cong K_
{ -1/2 }$. We cannot
prove however that the closure of the graph of this rational map
induces an inclusion $V\subset H^{k-1}(B,\QQ)\otimes H^1(A_K)$. In any
case, here is the construction.

Recall that $Y_k$ is the cyclic $3:1$ cover of $\PP^k$ totally
branched over a cubic hypersurface $X_{k-1}$. The main idea is that if
$l\subset \PP^{k}$ is a line, then its inverse image in $\PP^k$ is an
elliptic curve with an automorphism of order three with three fixed
points, so we can take that curve as $A_K$. Considering all lines
through a point (which we choose to be on $X_{k-1}$), we obtain an
isotrivial fibration on $Y_k$ which we trivialize after a 6:1 base
change; the base will be $B$.

In terms of formulae we have the following. Choose coordinates on
$\PP^k$ such that $(0:0:\ldots:0:1)\in X_k$ and write the equation of
$X_{k-1}$ as:
\[
X_{k-1}:\quad L x_k^2 + 2Q x_k + R=0,
\]
with $L,\,Q,\,R\in \CC[x_0,\ldots,x_{k-1}]$ homogeneous of degrees
$1,\,2,\,3$ respectively. Then $Y_k$ has equation
\[
Y_k:\quad x_{k+1}^3+L x_k^2 + 2Q x_k + R=0.
\]
Let $\bar{B}$ be the (singular) sextic subvariety of a $\PP^{k+1}$, with coordinates $x_0,\ldots,x_{k-1},\,x_{k+1},y$ defined by
\[
\bar{B}:\quad y^6+L^2(-LR+Q^2)=0,
\] 
and let $A_K$ be the elliptic curve defined by $u^3+v^2+1=0$. Then
we have a dominant rational map 
\[
\bar{B}\times A_K\longrightarrow Y_k,
\]
\[
\Bigl(\!(x_0\!:\ldots:x_{k-1}\!:\!x_{k+1}\!:\!y),\!(u,v)\!\Bigr)\!\mapsto\!
(x_0\!:\ldots:\!x_{k-1}\!:\!x_k\!:\!x_{k+1})\!=\!
(x_0\!:\ldots:\!x_{k-1}\!:{vy^3-LQ\over L^2}:{uy^2\over L}).
\]

\section{Torelli theorems via the half twist.}

\subsection{}
Let $X_{k-1}\subset \PP^k$ be a hypersurface such that $V\subset
H_0^k(Y_k,\ZZ)$ allows a half twist. Then we can associate to
$X_{k-1}$ two polarized Hodge structures of weight $k-1$, one would be
its primitive cohomology $H^{k-1}_0(X_{k-1},\QQ)$ and the other
$V_{1/2}$. We would like to know if one can recover $X_{k-1}$ from
$V_{1/2}$. We know this is the case for example if $k=2$ and $d=4$ by
the work of Kondo \cite{kondo99} and if $k=4$ and $d=3$ by recent work
of Allcock, Carlson and Toledo.

In case we do not only have the polarized Hodge structure $V_{1/2}$
but we also know the action of $\zeta = e^{ 2\pi i/d}$ on its
underlying $\ZZ$-module (this action comes from the action of $\alpha$
on $V$, which has the same underlying $\ZZ$-module), we show that $V$
can be recovered from $V_{1/2}$. In case $d$ is prime we thus recover
$H^k_0(Y_k,\ZZ)$ which in many case determines $Y_k$. In case $d$ is
not prime we do not know to what extent $V$ determines
$H^k_0(Y_k,\ZZ)$ and thus, via the Torelli theorems for hypersurfaces,
the variety $Y_k$.

\subsection{Lemma.}
Let $V$ be a Hodge structure of weight $k$ which allows a half twist $V_{1/2}$.
Then we recover $V$ from $V_{ 1/2 }$ in the following way
\[
V=(V_{1/2})_{-1/2}=\{ w \in V_{1/2} \otimes K_{-1/2} : \;
(\alpha\otimes\zeta^{-1})w=w \; \}.
\]

\ts
This follows from \ref{halftens} (cf.\ \cite{geemen01}). \qed

\subsection{Cubics.}
For the cubic $(k-1)$-fold $X_{k-1}$, the Hodge structures
$V_{1/2}(-1)$ and $W$ are equal because they have the same dimension
(see \ref{lemdim}). In this case we obtain a result which does not
require the knowledge of the automorphism of $V_{1/2}$.

\subsection{Lemma.}
\label{genitor} Suppose $k > 3$. The differential of the period map
which to $X_{ k-1 }$ associates the polarized Hodge structure $W$ is
  generically injective. In particular, this period map is generically
  finite.

\ts The tangent space to the space of periods for $W$ is the space of
compatible $(k-1)$-tuples of homomorphisms $W^{ p,r }\ra W^{ p-1,r+1
  }$ which preserve the polarization. The space of infinitesimal
  deformations of $X_{ k-1 }$ is isomorphic to
\[
H^0 \left(\frac{\cO_{\bP^k }}{ J'} (3)\right).
\]
Given a cubic polynomial $G_3$ (modulo $J'$), the image of the corresponding
deformation in the tangent space to the space of periods for $W$ is
the $(k-1)$-tuple of homomorphisms given by multiplication by $G_3$
in the description
\[
W^{ k+1 -p,p }\cong \frac{H^0 \left(\frac{\cO_{\bP^{ k+2 }}}{ J' } (-k + 3p
+3)\right)}{\left( x_{ k+1 } + x_{ k+2 } \right) H^0 (\cO_{\bP^k} (-k +3p
+2))}.
\]
To show that the differential is injective, it suffices to show that
one such homomorphism is non-zero, i.e., there is a $p$ and an element
of the above space which is not sent to zero by multiplication by $G_3$.
Again, as before, it is enough to check this for the Fermat cubic.
In this case, a basis of the above vector space is given by the
square-free monomials of degree $-k +3p +3$ in $x_0 ,\ldots ,x_{ k-1 }$
and products $x_{ k+1 }x_{ k+2 } M$ where $M$ is a square-free
monomial of degree $-k +3p +1$ in $x_0 ,\ldots x_{ k-1}$. An
infinitesimal deformation of the Fermat cubic can be uniquely
represented by a square-free cubic polynomial $G_3$ in the variables $x_0
,\ldots ,x_{ k-1 }$. After renaming the variables and possibly dividing
$G_3$ by a non-zero constant, we can assume that
\[
G_3 = x_0x_1x_2 + H_3
\]
where $H_3$ is a square-free cubic polynomial in the variables $x_0
,\ldots , x_{ k-1 }$ where the monomial $x_0x_1x_2$ does not appear.
Then, the image of any monomial $x_{ k+1 } x_{ k+2 } M$ where $M$ is a
square-free monomial in the variables $x_3 ,\ldots ,x_{ k-1}$ by
multiplication by $G_3$ is non-zero in $W^{ k-1 -p, p}$.
\qed

\providecommand{\bysame}{\leavevmode\hbox to3em{\hrulefill}\thinspace}

\end{document}